\documentclass{amsart}

\usepackage{amssymb}

\newtheorem{theorem}{Theorem}[section]
\newtheorem{corollary}[theorem]{Corollary}
\newtheorem{lemma}[theorem]{Lemma}

\newtheorem{example}[theorem]{Example}

\numberwithin{equation}{section}

\begin{document}

\title{Semiconjugacy of Quasiperiodic Flows and \\ Finite Index Subgroups of Multiplier Groups}


\email{\rm bakker@math.byu.edu}


\maketitle

\markboth{\sc L.F. Bakker}{\sc Semiconjugacy and Finite Index Subgroups}

\centerline{\scshape Lennard F. Bakker}
\medskip

{\footnotesize
\centerline{ Department of Mathematics, Brigham Young University}
\centerline{ Provo, UT 84602 USA }
}

\bigskip
\begin{quote}{\normalfont\fontsize{8}{10}\selectfont
{\bfseries Abstract.}
It will be shown that if $\phi$ is a quasiperiodic flow on the $n$-torus that is algebraic, if $\psi$ is a flow on the $n$-torus that is smoothly conjugate to a flow generated by a constant vector field, and if $\phi$ is smoothly semiconjugate to $\psi$, then $\psi$ is a quasiperiodic flow that is algebraic, and the multiplier group of $\psi$ is a finite index subgroup of the multiplier group of $\phi$. This will partially establish a conjecture that asserts that a quasiperiodic flow on the $n$-torus is algebraic if and only if its multiplier group is a finite index subgroup of the group of units of the ring of integers in a real algebraic number field of degree $n$.
\par}
\end{quote}

\section{Introduction}

The multiplier group of a (smooth, i.e.\ $C^\infty$) flow $\phi$ on the $n$-torus, $T^n={\mathbb R}^n/{\mathbb Z}^n$, is the image of the multiplier representation $\rho_\phi:S_\phi\to{\mathbb R}^*\equiv{\mathbb R}\setminus\{0\}$ of the flow's generalized symmetry group, $S_\phi = \{ R\in{\rm Diff}(T^n):{\rm there\ is\ }\alpha\in{\mathbb R}^*{\rm\ for\ which\ }R_*X_\phi=\alpha X_\phi \}$, which representation takes each $R\in S_\phi$ to its unique multiplier, or uniform scaling factor, $\alpha=\rho_\phi(R)$. (Here and elsewhere, ${\rm Diff}(T^n)$ is the group of smooth diffeomorphisms of $T^n$, $X_\phi(\theta)=(d/dt)\phi(t,\theta)\vert_{t=0}$, $\theta\in T^n$, is the smooth vector field generating $\phi$, $V_*X_\phi={\bf T}V X_\phi V^{-1}$ is the push-forward of $X_\phi$ by $V\in{\rm Diff}(T^n)$, and ${\bf T}$ is the tangent functor.) The multiplier group, $M_\phi=\rho_\phi(S_\phi)$, is an absolute (but incomplete) invariant of the smooth conjugacy class of $\phi$ (Theorem 4.2 in \cite{BC}), absolute in the sense that if $\phi$ is smoothly conjugate to $\psi$ (i.e.\ if there exists a $V\in{\rm Diff}(T^n)$ such that $V_*X_\phi=X_\psi$), then $M_\phi=M_\psi$ (they are identical subsets of ${\mathbb R}^*$, the multiplicative real group). It will be shown that if $\phi$ is a quasiperiodic flow on $T^n$, if $\psi$ is a flow on $T^n$ that is smoothly conjugate to a flow generated by a constant vector field, and if $\phi$ is smoothly semiconjugate to $\psi$ (i.e.\ if there exists a surjective $V$ in $C^\infty(T^n)$ -- the monoid of smooth maps from $T^n$ to $T^n$ -- such that ${\bf T}V X_\phi=X_\psi V$),  then
\begin{enumerate}
\item $\psi$ is quasiperiodic (Theorem \ref{covering}),
\item $M_\phi\supset M_\psi$ (Corollary \ref{subgroup}), and
\item $M_\psi$ is a finite index subgroup of $M_\phi$ whenever $M_\phi$ is finitely generated (Corollary \ref{finite}).
\end{enumerate}
$\big($Recall that a flow $\phi$ on $T^n$ is quasiperiodic if and only if there is $V\in{\rm Diff}(T^n)$ such that $V_*X_\phi$ is a constant vector field whose $n$ real components, or frequencies, are independent over ${\mathbb Q}$ (pp.79-80 in \cite{HB})$\big)$. The {\it main result} (Theorem \ref{algebraic}) is that if $\phi$ is a quasiperiodic flow on $T^n$ that is algebraic (i.e.\ there exists a $\vartheta\in{\mathbb R}^*$ and a real algebraic number field $F$ of degree $n$ such that the components of the constant vector field $\vartheta(V_*X_\phi)$ form a basis for $F$ as an $n$-dimensional vector space over ${\mathbb Q}$), if $\psi$ is a flow that is smoothly conjugate to a flow generated by a constant vector field, and if $\phi$ is smoothly semiconjugate to $\psi$, then $\psi$ is a quasiperiodic flow that is algebraic, and $M_\psi$ is a finite index subgroup of $M_\phi$. (Recall that an algebraic number field is a finite dimensional field extension of ${\mathbb Q}$.)

\begin{example}\label{example}{\rm (Illustration of Main Theorem) Consider the flows $\phi$ and $\psi$ on $T^2$, equipped with global coordinates $\theta=(\theta_1,\theta_2)$, generated by
\[ X_\phi=\frac{\partial}{\partial \theta_1}+\big(1+\sqrt 2\big) \frac{\partial}{\partial \theta_2}{\rm\ \ and\ \ } X_\psi=\big(4+\sqrt 2\big)\frac{\partial}{\partial \theta_1}+\big(3+2\sqrt 2\big)\frac{\partial}{\partial \theta_2}.\]
The flow $\phi$ is quasiperiodic (because $(1+\sqrt 2)/1$, the ratio of its frequencies, is a root of the irreducible quadratic polynomial $z^2-2z-1$ in the polynomial ring ${\mathbb Q}[z]$; see Theorem 2.3 in \cite{BA1}). Also, the components of $X_\phi$ form a basis for the real algebraic number field $F={\mathbb Q}(\sqrt 2)=\{a+b\sqrt 2:a,b\in{\mathbb Q}\}$ of degree $2$; hence $\phi$ is algebraic, and $M_\phi\subset {\mathfrak o}^*_F$ (Theorem 3.4 in \cite{BA3}). Let $V\in C^\infty(T^n)$ be the surjective map induced by the invertible integer matrix
\[ \begin{bmatrix} 3 & 1 \\ 1 & 2\end{bmatrix}.\]
Then $\phi$ is semiconjugate to $\psi$ by $V$ because ${\bf T}V X_\phi=X_\psi=X_\psi V$. The flow $\psi$ is also quasiperiodic (because the ratio $(3+2\sqrt 2)/(4+\sqrt 2)=(8+5\sqrt 2)/14$ is a root of the irreducible $z^2-(8/7)z+(1/14)$ in ${\mathbb Q}[z]$), and the components of $X_\psi$ form a basis for ${\mathbb Q}(\sqrt 2)$ as well; hence $\psi$ is algebraic and $M_\psi\subset{\mathfrak o}^*_F$ too.  The unimodular matrix
\[ \begin{bmatrix}0 & 1\\1 & 2\end{bmatrix}\]
induces an $R\in{\rm Diff}(T^2)$ which satisfies $R_* X_\phi=\big(1+\sqrt 2\big)X_\phi$, so that $R\in S_\phi$ with $\rho_\phi(R)=1+\sqrt 2\in M_\phi$. It can be shown that $M_\phi=\{\pm(1+\sqrt 2)^k:k\in{\mathbb Z}\}={\mathfrak o}^*_F$ (see Example 3.9 in \cite{BA3}), so that $M_\psi\subset M_\phi$ with $M_\phi$ finitely generated. Since for each $\alpha\in M_\psi$ there are $u_1,u_2\in{\mathbb Z}$ such that
\[ \alpha = u_1\bigg(\frac{4+\sqrt 2}{3+2\sqrt 2}\bigg) + u_2 = u_1(8-5\sqrt 2) + u_2\]
(Corollary 4.5 in \cite{BA2}), it follows that neither $1+\sqrt 2$ nor $(1+\sqrt 2)^2=3+2\sqrt 2$ belong to $M_\psi$. On the other hand, the unimodular matrix
\[ \begin{bmatrix}-1 & 14 \\ -1 & 15\end{bmatrix} \]
induces a $Q\in{\rm Diff}(T^2)$ which satisfies $Q_*X_\psi=(1+\sqrt 2)^3 X_\psi$, so that $Q\in S_\psi$ and $(1+\sqrt 2)^3=7+5\sqrt 2=\rho_\psi(Q)\in M_\psi$. (Notice that $V$, $R$, and $Q$ satisfies
\[ {\bf T}V({\bf T}R)^3=\begin{bmatrix} 3 & 1 \\ 1 & 2\end{bmatrix}\begin{bmatrix} 0 & 1 \\ 1 & 2\end{bmatrix}^3=\begin{bmatrix}-1 & 14 \\ -1 & 15\end{bmatrix}\begin{bmatrix} 3 & 1 \\ 1 & 2\end{bmatrix}={\bf T}Q{\bf T}V,\]
so that $({\bf T}R)^3$ is ``semiconjugate'' to ${\bf T}Q$ by ${\bf T}V$.) Since $\pm 1\in M_\psi$ (which is the case for any quasiperiodic flow by Theorem 2.3 in \cite{BA2}), then $M_\psi=\{\pm(1+\sqrt 2)^{3k}:k\in{\mathbb Z}\}$ which is a subgroup of $M_\phi$ of index $3$.
}\end{example}

As in the example there is, for a given real algebraic number field $F$ of degree $n$, a quasiperiodic flow $\phi$ on $T^n$ such that $\phi$ is algebraic and $M_\phi={\mathfrak o}^*_F$ (the frequencies of $\phi$ are the elements of an integral basis for $F$; see Theorem 3.8 in \cite{BA3}). By Dirichlet's Unit Theorem (see p.21 in \cite{SD}), the group ${\mathfrak o}^*_F$ is a finitely generated (abelian) group that contains at least one infinite cyclic factor (so that $M_\phi\setminus\{1,-1\}\ne\emptyset$, 
i.e.\ $\phi$ possesses multipliers other than $\pm 1$). Applying the main result to any flow $\psi$ on $T^n$ which is smoothly conjugate to a flow generated by a constant vector field,  for which $\phi$ is smoothly semiconjugate to $\psi$, shows that $\psi$ is a quasiperiodic flow that is algebraic and that $M_\psi$ is a finite index subgroup of ${\mathfrak o}^*_F$. This partially establishes a conjecture (Conjecture 4.2 in \cite{BA3}) that asserts that a quasiperiodic flow $\psi$ on $T^n$ is algebraic if and only if there is a real algebraic number field $F$ of degree $n$ such that $M_\psi$ is a finite index subgroup of ${\mathfrak o}^*_F$.

\section{Semiconjugacy and Covering Maps}

One of the key elements in the proof of main result is the topological result that a surjective $V\in C^\infty(T^n)$ is a covering map (Theorem \ref{covering}) whenever $V$ is a smooth semiconjugacy from $\phi$ to $\psi$ (i.e.\ ${\bf T}VX_\phi=X_\psi V$) with $\phi$ quasiperiodic and $\phi$ smoothly semiconjugate to a flow generated by a constant vector field. (Recall that a $V\in C^\infty(T^n)$ is a covering map if each point in $T^n$ has a connected open neighbourhood $O$ whose inverse image, $V^{-1}(O)$, is a disjoint union of open sets in $T^n$ each homeomorphic to $O$ by $V$.) 

The (universal) covering map $\pi:{\mathbb R}^n\to T^n$ is a local diffeomorphism whose deck transformations are the translations in ${\mathbb R}^n$ by $m\in{\mathbb Z}^n$. A lift of a continuous map $K:T^n\to T^n$ is a continuous map $\hat K:{\mathbb R}^n\to{\mathbb R}^n$ such that $K\pi=\pi\hat K$. (Such a $\hat K$ exists by the Lifting Theorem (see p.143 in \cite{BR} for example) since the fundamental group of ${\mathbb R}^n$ is trivial; any two lifts of $K$ differ by a deck transformation of $\pi$.) Any lift of $K\in C^\infty(T^n)$ is in $C^\infty({\mathbb R}^n)$ -- the monoid of smooth maps from ${\mathbb R}^n$ to ${\mathbb R}^n$ -- since $K$ is smooth and $\pi$ is a local diffeomorphism. The identification of the $n$-torus with ${\mathbb R}^n/{\mathbb Z}^n$ determines the global coordinates $\theta=(\theta_1,\dots,\theta_n)$ on $T^n$, and the concomitant global coordinates $x=(x_1,x_2,\dots,x_n)$ on ${\mathbb R}^n$ such that $\theta_i=x_i$ mod $1$, $i=1,2,\dots,n$. In terms of these coordinates on ${\mathbb R}^n$, let ${\rm M}(n,{\mathbb Z})$ denote the monoid of $n\times n$ matrices with integer entries, and let ${\rm GL}(n,{\mathbb Z})$ denote the unimodular group (i.e.\ those matrices in ${\rm M}(n,{\mathbb Z})$ with determinant $\pm 1$). The proof of the following relationship between a map in $C^\infty(T^n)$ and its lifts in $C^\infty({\mathbb R}^n)$ uses standard arguments in topology, and is therefore omitted. (Here, ${\rm Diff}({\mathbb R}^n)$ is the group of diffeomorphisms of ${\mathbb R}^n$.)

\begin{theorem}\label{liftmap} A map $W$ in $C^\infty({\mathbb R}^n)$ $\big(in {\rm\ Diff}({\mathbb R}^n)\big)$ is a lift of a map in $C^\infty(T^n)$ $\big(in {\rm\ Diff}(T^n)\big)$ if and only if there exists $B$ in ${\rm M}(n,{\mathbb Z})$ $\big(in {\rm\ GL}(n,{\mathbb Z})\big)$ such that $W(x+m)-W(x)=Bm$ for all $x\in{\mathbb R}^n$ and all $m\in{\mathbb Z}^n$.
\end{theorem}

The covering map $\pi:{\mathbb R}^n\to T^n$ induces the covering map ${\bf T}\pi:{\bf T}{\mathbb R}^n\to {\bf T}T^n$. There is only one lift of a vector field $X$ on $T^n$ that is a vector field on ${\mathbb R}^n$ (Lemma 3.2 in \cite{BA2}). This lift, denoted by $\hat X$, satisfies $X\pi={\bf T}\pi \hat X$. In particular, the lift of a constant vector field $X=\sum_{i=1}^n a_i\partial/\partial\theta_i$, $a_i\in{\mathbb R}$, on $T^n$ is the constant vector field $\hat X=\sum_{i=1}^n a_i\partial/\partial x_i$ on ${\mathbb R}^n$. The restriction that the equation ${\bf T}V X_\phi=X_\psi V$ on ${\bf T}T^n$ places on $V$ is extracted by first lifting this equation to ${\bf T}{\mathbb R}^n$.

\begin{lemma}\label{charliftpde} Let $X$ and $Y$ be vector fields on $T^n$, and let $\hat K\in C^\infty({\mathbb R}^n)$ be a lift of $K\in C^\infty(T^n)$. Then ${\bf T}KX=YK$ if and only if ${\bf T}\hat K\hat X=\hat Y\hat K$.
\end{lemma}

\begin{proof} Since $K\pi=\pi\hat K$, $X\pi={\bf T}\pi \hat X$, and $Y\pi={\bf T}\pi \hat Y$, it follows that
\[{\bf T}\pi\big({\bf T}\hat K \hat X-\hat Y\hat K) = \big( {\bf T}KX-YK\big)\pi.\]
If ${\bf T}KX=YK$, then ${\bf T}\pi({\bf T}\hat K\hat X-\hat Y\hat K)=0$, so that ${\bf T}\hat K\hat X=\hat Y\hat K$ since $\pi$ is a local diffeomorphism. If ${\bf T}\hat K\hat X=\hat Y\hat K$, then $({\bf T}KX-YK)\pi=0$, so that ${\bf T}KX=YK$ since $\pi$ is surjective. 
\end{proof}

\begin{theorem}\label{solvepde} Let $X$ and $Y$ be constant vector fields on $T^n$, and let $\hat K$ in $C^\infty({\mathbb R}^n)$ $\big(in\ {\rm Diff}({\mathbb R}^n \big)$ be a lift of $K$ in $C^\infty(T^n)$ $\big(in\ {\rm Diff}(T^n)\big)$. If the components of $X$ are independent over ${\mathbb Q}$ and if ${\bf T}K X=Y$, then there exist $B$ in ${\rm M}(n,{\mathbb Z})$ $\big(in\ {\rm GL}(n,{\mathbb Z})\big)$ and $c\in{\mathbb R}^n$ such that $\hat K(x)=Bx+c$ for all $x\in{\mathbb R}^n$.
\end{theorem}

\begin{proof} Suppose that ${\bf T}K X=Y$ where $X = \sum_{i=1}^n a_i\partial/\partial \theta_i$ and $Y = \sum_{i=1}^n d_i \partial/\partial \theta_i$ for constants $a_i,d_i\in{\mathbb R}$. Then ${\bf T}K X=YK$ so that ${\bf T}\hat K\hat X=\hat Y\hat K$ by Lemma \ref{charliftpde}; hence ${\bf T}\hat K\hat X=\hat Y K=\hat Y$ since $\hat Y$ is a constant vector field. The form of $\hat K$ is
\[ \hat K(x_1,x_2,\dots,x_n)=\big(f_1(x_1,x_2,\dots,x_n),\dots,f_n(x_1,x_2,\dots,x_n)\big)\]
for smooth functions $f_i:{\mathbb R}^n\to{\mathbb R}$, $i=1,2,\dots,n$. The equation ${\bf T}\hat K\hat X=\hat Y$ written out is an uncoupled system of $n$ linear, first order partial differential equations with constant coefficients:
\[ \sum_{j=1}^n a_j\frac{\partial f_i}{\partial x_j} = d_i, \ \ i=1,\dots,n.\]

Suppose the coefficients of $X$ are linearly independent over ${\mathbb Q}$. This implies that $a_n\ne0$, so that by the method of characteristics (see \cite{FJ} for example), the general solution of the uncoupled system is
\[ f_i(x)=\frac{d_i}{a_n}x_n+ h_i\bigg(x_1-\frac{a_1}{a_n}x_n,\dots,x_{n-1}-\frac{a_{n-1}}{a_n}x_n\bigg)\]
for arbitrary smooth functions $h_i:{\mathbb R}^{n-1}\to{\mathbb R}$, $i=1,\dots,n$.
Since $\hat K$ in $C^\infty({\mathbb R}^n)$ $\big({\rm in\ Diff}({\mathbb R}^n)\big)$ is a lift of $K$ in $C^\infty(T^n)$ $\big({\rm in\ Diff}(T^n)\big)$, there is by Theorem \ref{liftmap}, a matrix $B=(b_{ij})$ in ${\rm M}(n,{\mathbb Z})$ $\big({\rm in\ GL}(n,{\mathbb Z})\big)$ such that 
\[ f_i(x+m)-f_i(x)=\sum_{j=1}^n b_{ij}m_j\]
for all $x\in{\mathbb R}^n$ and all $m=(m_1,\dots,m_n)\in{\mathbb Z}^n$. Set $s_j=x_j-(a_j/a_n)x_n$ for $j=1,\dots,n-1$. Then each $h_i$ satisfies
\begin{align*} \sum_{j=1}^n b_{ij}m_j 
& = h_i\bigg(s_1+m_1-\frac{a_1}{a_n}m_n,\dots,s_{n-1}+m_{n-1}-\frac{a_{n-1}}{a_n}m_n\bigg) \\ 
& \ \ \ \ -h_i(s_1,\dots,s_{n-1}) + \frac{d_i}{a_n}m_n
\end{align*}
for all $(s_1,\dots,s_{n-1})\in{\mathbb R}^{n-1}$ and all $m\in{\mathbb Z}^n$. Differentiation of both sides of  this with respect to $s_j$ followed by evaluation at $s=0$ gives for each $i=1,\dots,n$ and each $j=1,\dots,n-1$ that
\[ \frac{\partial h_i}{\partial s_j}\bigg(m_1-\frac{a_1}{a_n}m_n,\dots, m_{n-1}-\frac{a_{n-1}}{a_n}m_n\bigg)=\frac{\partial h_i}{\partial s_j}(0,\dots,0)\]
for all $m\in{\mathbb Z}^n$. The independence of the $a_i$'s over ${\mathbb Q}$ implies the denseness of subset $\big\{\big(m_1-(a_1/a_n)m_n,\dots,m_{n-1}-(a_{n-1}/a_n)m_n\big):m\in{\mathbb Z}^n \big\}$ in ${\mathbb R}^{n-1}$ (by Lemma 4.2 in \cite{BA2}), so that the smoothness of the $h_i$'s implies that $\partial h_i/\partial s_j=b^\prime_{ij}$ are real constants for $i=1,\dots,n$, $j=1,\dots,n-1$. By Taylor's Theorem, there are constants $c_i\in{\mathbb R}$ such that
\[ h_i(s_1,\dots,s_{n-1})=c_i+\sum_{j=1}^{n-1}b^\prime_{ij}s_j\]
for each $i=1,\dots,n$. Set $b^\prime_{in}=d_i/a_n-\sum_{j=1}^{n-1}b^\prime_{ij}(a_j/a_n)$. Then
\[ f_i(x_1,\dots,x_n)= c_i + \frac{d_i}{a_n}x_n+\sum_{j=1}^{n-1} b^\prime_{ij}\bigg(x_j-\frac{a_j}{a_n}x_n\bigg) = c_i +\sum_{j=1}^n b^\prime_{ij}x_j\]
for each $i=1,\dots,n$. Since $f_i(x+m)-f(x)=\sum_{j=1}^n b_{ij}m_j$ for all $x\in{\mathbb R}^n$ and all $m\in{\mathbb Z}^n$, then $b^\prime_{ij}=b_{ij}$ for all $i,j=1,\dots,n$, so that $\hat K(x)=Bx+c$ with $c=(c_1,\dots,c_n)\in{\mathbb R}^n$.
\end{proof}

An immediate consequence of Theorem \ref{solvepde} is that up to a translation, a $V\in C^\infty(T^n)$ (not assumed surjective) that satisfies ${\bf T}VX_\phi=X_\psi$ for $X_\phi$ and $X_\psi$ constant vector fields and $\phi$ quasiperiodic, is a map that is induced by a matrix $B$ in ${\rm M}(n,{\mathbb Z})$, i.e.\ ${\bf T}\hat V=B$. Those $B$ in ${\rm M}(n,{\mathbb Z})$ with $\det B\ne 0$ belong to ${\rm GL}(n,{\mathbb Q})$ (the group of $n\times n$ matrices with rational entries and nonzero determinant).

\begin{lemma}\label{tangentofV} Let $X$ and $Y$ be constant vector fields on $T^n$, and let $\hat V$ be the lift of a surjective $V\in C^\infty(T^n)$. If the components of $X$ are independent over ${\mathbb Q}$ and if ${\bf T}VX=Y$, then ${\bf T}\hat V \in{\rm M}(n,{\mathbb Z})\cap{\rm GL}(n,{\mathbb Q})$.
\end{lemma}

\begin{proof} Suppose that $V$ is surjective and that ${\bf T}VX =Y$. The independence of the components of $X$ over ${\mathbb Q}$ implies by Theorem \ref{solvepde} that $\hat V(x)=Bx+c$ for $B\in{\rm M}(n,{\mathbb Z})$ and $c\in{\mathbb R}^n$. Thus $B={\bf T}\hat V\in {\rm M}(n,{\mathbb Z})$. Suppose that $\det B=0$. Then $\hat V({\mathbb R}^n)$ is a nowhere dense subset of ${\mathbb R}^n$. If $\hat V({\mathbb R}^n)+m$ is the translation of $\hat V({\mathbb R}^n)$ by $m\in{\mathbb Z}^n$, then ${\mathbb R}^n\ne\cup_{m\in{\mathbb Z}^n} \big(\hat V({\mathbb R}^n)+m\big)$ by the Baire Category Theorem. This means that there exists $x\in{\mathbb R}^n$ such that $x+m\not\in \hat V({\mathbb R}^n)$ for all $m\in{\mathbb Z}^n$. Thus $\pi(x)\not\in \pi\hat V({\mathbb R}^n)=V\pi({\mathbb R}^n)=V(T^n)$, so that $V$ is not surjective. This contradiction implies that ${\bf T}\hat V\in{\rm GL}(n,{\mathbb Q})$.
\end{proof}

\begin{theorem}\label{covering} If $\phi$ is a quasiperiodic flow on $T^n$, if $\psi$ is a flow on $T^n$ that is smoothly conjugate to a flow generated by a constant vector field, and if $\phi$ is smoothly semiconjugate to $\psi$ by a surjective $V\in C^\infty(T^n)$, then $\psi$ is quasiperiodic and $V$ is a covering map. If, in addition, $\phi$ is algebraic, then $\psi$ is algebraic also.
\end{theorem}

\begin{proof} By the hypotheses, there are $K,U\in{\rm Diff}(T^n)$ and a surjective $V\in C^\infty(T^n)$ such that $K_*X_\phi$ and $U_*X_\psi$ are constant vector fields, the components of $X_\phi$ being independent over ${\mathbb Q}$, and ${\bf T}VX_\phi=X_\psi V$.

Let $\Phi$ be the quasiperiodic flow on $T^n$ generated by $X_\Phi\equiv K_*X_\phi$, and let $\Psi$ be the quasiperiodic flow on $T^n$ generated by $X_\Psi\equiv U_*X_\psi$. Then $\Phi$ is smoothly semiconjugate to $\Psi$ by the surjective $UVK^{-1}\in C^\infty(T^n)$ because ${\bf T}(UVK^{-1})X_\Phi={\bf T}U{\bf T}V X_\phi K^{-1}={\bf T}UX_\psi VK^{-1}=X_\Psi(UVK^{-1})=X_\Psi$. Since a lift of $UVK^{-1}$ is $\hat U\hat V\hat K^{-1}$, it follows by Lemma \ref{tangentofV} that $B={\bf T}(\hat U\hat V\hat K^{-1})\in {\rm M}(n,{\mathbb Z})\cap{\rm GL}(n,{\mathbb Q})$. Since $X_\Psi=B X_\Phi$ by Lemma \ref{charliftpde}, the independence of the components of $X_\Phi$ over ${\mathbb Q}$ then implies that of the components of $X_\Psi$; thus $\psi$ is quasiperiodic.

If $\phi$ is algebraic, then there is a $\vartheta\in{\mathbb R}^*$ and a real algebraic number field $F$ of degree $n$ such that the components of $\vartheta X_\Phi$ form a basis for $F$. Since $B\in{\rm GL}(n,{\mathbb Q})$ and since $X_\Psi=B X_\Phi$, the components of $\vartheta X_\Psi$ form a basis for $F$ as well. Hence $\psi$ is algebraic.

Since $UVK^{-1}$ is induced by $B$ (up to a translation) and since $\det B\in{\mathbb Z}\setminus\{0\}$, then $UVK^{-1}$ is a $\vert\det B\vert$ to $1$ map. This means that for any sufficiently small open connected subset $O$ of $T^n $, the inverse image $(UVK^{-1})^{-1}(O)$ consists of $\vert\det B\vert$ components each homeomorphic to $O$ by $UVK^{-1}$, and hence $UVK^{-1}$ is a covering map. Since $U,K\in {\rm Diff}(T^n)$, it follows that $V$ is a covering map as well.
\end{proof}

\section{Semiconjugacy and Multiplier Groups}

Additional key elements in the proof of the main result are two ``semiconjugacy'' relations between $S_\phi$ and $S_\psi$ in terms of a smooth semiconjugacy $V$ from $\phi$ to $\psi$ (described in Theorems \ref{semiconjugate} and \ref{power}), which are valid when $\phi$ and $\psi$ are generated by constant vector fields and $\phi$ is quasiperiodic.

\begin{theorem}\label{semiconjugate} Let $\phi$ and $\psi$ be flows on $T^n$ generated by  constant vector fields, with $\phi$ quasiperiodic. If $\phi$ is smoothly semiconjugate to $\psi$ by a surjective $V\in C^\infty(T^n)$, then for each $Q\in S_\psi$ there is $R\in S_\phi$ such that $QV=VR$ and $\rho_\phi(R)=\rho_\psi(Q)$.
\end{theorem}

\begin{proof} Suppose ${\bf T}V X_\phi=X_\psi$ for a surjective $V\in C^\infty(T^n)$. By Theorem \ref{covering}, the quasiperiodicity of $\phi$ implies that $\psi$ is quasiperiodic and that $V$ is a covering map. A fixed lift $\hat V$ satisfies ${\bf T}\hat V\hat X_\phi=\hat X_\psi$ by Lemma \ref{charliftpde} where ${\bf T}\hat V\in{\rm M}(n,{\mathbb Z})\cap{\rm GL}(n,{\mathbb Q})$ by Lemma \ref{tangentofV}. Since ${\bf T}V{\bf T}\pi={\bf T}\pi{\bf T}\hat V$, it follows that $V$ is a local diffeomorphism.

Let $Q\in S_\psi$ and set $\alpha=\rho_\psi(Q)$. Since $V$ is a covering map, there is by the Lifting Theorem (see p.143 in \cite{BR} for example) a continuous map $R:T^n\to T^n$ such that $QV=VR$. Since $Q\in {\rm Diff}(T^n)$ and since $V$ is a local diffeomorphism, it follows that $R\in C^\infty(T^n)$. Fix a lift $\hat R$. If $\hat Q$ is a lift of $Q$, then
\[ \pi \hat Q\hat V=QV\pi=VR\pi=\pi \hat V\hat R,\]
so that $\hat Q\hat V$ and $\hat V\hat R$ are lifts of $QV$. These two lifts differ by a deck transformation: there is an $m\in{\mathbb Z}^m$ such that $i_m\hat Q\hat V=\hat V\hat R$ where $i_m(x)=x+m$. But $i_m\hat Q$ is just another lift of $Q$. Replacing $\hat Q$ with $i_m\hat Q$ gives $\hat Q\hat V=\hat V\hat R$ and ${\bf T}(\hat Q\hat V)={\bf T}(\hat V\hat R)$.

The generalized symmetry $Q$ of $\psi$ satisfies ${\bf T}Q X_\psi=Q_*X_\psi=\alpha X_\psi$; hence ${\bf T}\hat Q \hat X_\psi = \alpha \hat X_\psi$ by Lemma \ref{charliftpde}. (The only lift of $\alpha X_\psi$ that is a vector field is $\alpha\hat X_\psi$; see Lemma 3.4 in \cite{BA2}.)  The quasiperiodicity of $\psi$ implies by Theorem \ref{solvepde} that ${\bf T}\hat Q\in {\rm GL}(n,{\mathbb Z})$. Since the map $\hat R$ satisfies
\begin{align*}
{\bf T}\hat V({\bf T}\hat R \hat X_\phi)
& = {\bf T}(\hat V\hat R)\hat X_\phi = {\bf T}(\hat Q\hat V)\hat X_\phi = {\bf T}\hat Q{\bf T}\hat V \hat X_\phi \\
& = {\bf T}\hat Q \hat X_\psi = \alpha \hat X_\psi = \alpha {\bf T}\hat V \hat X_\phi = {\bf T}\hat V(\alpha \hat X_\phi),
\end{align*}
and since $\det {\bf T}\hat V\ne 0$, it follows that ${\bf T}\hat R \hat X_\phi=\alpha \hat X_\phi$; hence ${\bf T}R X_\phi=\alpha X_\phi$ by Lemma \ref{charliftpde}. The quasiperiodicity of $\phi$ implies by Theorem \ref{solvepde} that ${\bf T}\hat R\in {\rm M}(n,{\mathbb Z})$. Since ${\bf T}(\hat Q\hat V)={\bf T}(\hat V\hat R)$, it follows that
\[ \det {\bf T}\hat Q \det {\bf T}\hat V=\det {\bf T}\hat V \det {\bf T}\hat R.\]
Since $\det {\bf T}\hat Q=\pm 1$ and $\det{\bf T}\hat V\ne0$, then $\det {\bf T}\hat R=\pm1$, so that ${\bf T}\hat R\in {\rm GL}(n,{\mathbb Z})$. With $B={\bf T}\hat R$, then $\hat R(x)=Bx+c$ for $c\in{\mathbb R}^n$, so that $\hat R(x+m)-\hat R(x)=Bm$ for all $x\in{\mathbb R}^n$ and all $m\in{\mathbb Z}^n$. Thus $R\in {\rm Diff}(T^n)$ by Theorem \ref{liftmap}. Since $R_*X_\phi={\bf T}R X_\phi=\alpha X_\phi$, it follows that $R\in S_\phi$ with $\rho_\phi(R)=\alpha=\rho_\psi(Q)$.
\end{proof}

\begin{corollary}\label{subgroup} If $\phi$ is a quasiperiodic flow on $T^n$, if $\psi$ is a flow on $T^n$ that is smoothly conjugate to a flow generated by a constant vector field, and if $\phi$ is smoothly semiconjugate to $\psi$, then $M_\phi\supset M_\psi$.
\end{corollary}

\begin{proof} Since the multiplier group is an absolute invariant of the smooth conjugacy class of a flow, and since a quasiperiodic flow is smoothly conjugate to a flow generated by a constant vector field, Theorem \ref{semiconjugate} implies that $M_\phi\supset M_\psi$.
\end{proof}

\begin{theorem}\label{power} Let $\phi$ and $\psi$ be flows on $T^n$ generated by constant vector fields, with $\phi$ quasiperiodic. If $\phi$ is smoothly semiconjugate to $\psi$ by a surjective $V\in C^\infty(T^n)$, then for each $R\in S_\phi$ there exists $k\in{\mathbb Z}^+$ and $Q\in S_\psi$ such that $QV=VR^k$ and $\rho_\psi(Q)=[\rho_\phi(R)]^k$.
\end{theorem}

\begin{proof} Suppose ${\bf T}V X_\phi=X_\psi$ for a surjective $V\in C^\infty(T^n)$, and fix a lift $\hat V\in C^\infty({\mathbb R}^n)$. The quasiperiodicity of $\phi$ implies by Theorem \ref{solvepde} that there is $B\in{\rm M}(n,{\mathbb Z})$ and $c\in{\mathbb R}^n$ such that $\hat V(x)=Bx+c$. By Lemma \ref{tangentofV}, the matrix $B={\bf T}\hat V$ is in ${\rm GL}(n,{\mathbb Q})$, so that $\hat V\in{\rm Diff}(T^n)$. Since ${\bf T}\hat V\hat X_\phi=\hat X_\psi$ by Lemma \ref{charliftpde}, then ${\bf T}\hat V^{-1}\hat X_\psi=\hat X_\phi$.

Let $R\in S_\phi$, set $\alpha=\rho_\phi(R)$, and fix a lift $\hat R$. Since ${\bf T}R X_\phi=R_*X_\phi=\alpha X_\phi$, then ${\bf T}\hat R \hat X_\phi=\alpha\hat X_\phi$ by Lemma \ref{charliftpde}, and ${\bf T}\hat R\in{\rm GL}(n,{\mathbb Z})$ by the quasiperiodicity of $\phi$ and Theorem \ref{solvepde}. Since ${\bf T}\hat V\in{\rm GL}(n,{\mathbb Q})$ and since ${\rm GL}(n,{\mathbb Z})$ is a subgroup of ${\rm GL}(n,{\mathbb Q})$, the matrix ${\bf T}\hat V{\bf T}\hat R{\bf T}\hat V^{-1}$ is in ${\rm GL}(n,{\mathbb Q})$ and has the same characteristic polynomial as that of ${\bf T}\hat R$, which polynomial has integer coefficients. (The matrix ${\bf T}\hat V{\bf T}\hat R{\bf T}\hat V^{-1}$ would be in ${\rm GL}(n,{\mathbb Z})$ for any $R\in S_\phi$ if ${\rm GL}(n,{\mathbb Z})$ were a normal subgroup of ${\rm GL}(n,{\mathbb Q})$; but as simple examples show, ${\rm GL}(n,{\mathbb Z})$ is not a normal subgroup of ${\rm GL}(n,{\mathbb Q})$.) A result in commutative ring theory (Theorem 2 in \cite{DRTW}) states that an $n\times n$ matrix with rational entries, with determinant $\pm 1$, and with a characteristic polynomial having integer coefficients, has a power whose entries are integers. Thus, there is a $k\in{\mathbb Z}^+$ such that $\tilde B=({\bf T}\hat V{\bf T}\hat R{\bf T}\hat V^{-1})^k={\bf T}\hat V{\bf T}\hat R^k{\bf T}\hat V^{-1}\in{\rm GL}(n,{\mathbb Z})$. The map $W:{\mathbb R}^n\to{\mathbb R}^n$ defined by $W=\hat V\hat R^k\hat V^{-1}$ is in ${\rm Diff}({\mathbb R}^n)$ (since $\hat V$ and $\hat R$ are) and satisfies ${\bf T}W=\tilde B$. Thus there is $\tilde c\in{\mathbb R}^n$ such that $W(x)=\tilde Bx+\tilde c$, from which it follows that $W(x+m)-W(x)=\tilde Bm$ for all $x\in{\mathbb R}^n$ and all $m\in{\mathbb Z}^n$. By Theorem \ref{liftmap}, the map $W$ is a lift of a $Q\in{\rm Diff}(T^n)$. Thus $\hat Q=\hat V\hat R^k\hat V^{-1}$, so that
\[ QV\pi=\pi\hat Q\hat V=\pi\hat V\hat R^k=VR^k\pi.\]
The surjectiveness of $\pi$ implies that $QV=VR^k$. Moreover,
\[ {\bf T}\hat Q \hat X_\psi
= {\bf T}\hat V{\bf T}\hat R^k{\bf T}\hat V^{-1} \hat X_\psi = {\bf T}\hat V{\bf T}\hat R^k \hat X_\phi = {\bf T}\hat V(\alpha^k \hat X_\phi) = \alpha^k \hat X_\psi,\]
so that $Q_*X_\psi={\bf T}QX_\psi=\alpha^k X_\psi$ by Lemma \ref{charliftpde}. Therefore, $Q\in S_\psi$ and $\rho_\psi(Q)=\alpha^k=[\rho_\phi(R)]^k$.
\end{proof}

\begin{corollary}\label{finite} Suppose that $\phi$ is a quasiperiodic flow on $T^n$, that $\psi$ is a flow on $T^n$ that is smoothly conjugate to a flow generated by a constant vector field, and that $\phi$ is smoothly semiconjugate to $\psi$. If $M_\phi$ is finitely generated, then $M_\psi$ is a finite index subgroup of $M_\phi$.
\end{corollary}

\begin{proof} It suffices to show that $M_\psi$ is a finite index subgroup of $M_\phi$ when $\phi$ and $\psi$ are generated by constant vector fields (since the multiplier group is an absolute invariant of the smooth conjugacy class of a flow, and since a quasiperiodic flow is smoothly conjugate to a flow generated by a constant vector field). 

Suppose that $M_\phi$ is finitely generated: there are $\alpha_i\in M_\phi$, $i=1,\dots,l$ such that $M_\phi=\langle a_1,\dots,a_l\rangle$ (the smallest subgroup of ${\mathbb R}^*$ that contains $\{a_1,\dots,a_l\}$). For each $\alpha_i$ there is by Theorem \ref{power} a $k_i\in{\mathbb Z}^+$ such that $\alpha_i^{k_i}\in M_\psi$. Thus $G=\langle \alpha_1^{k_1},\dots,\alpha_l^{k_l}\rangle$ is a subgroup of $M_\psi$. Since $M_\psi$ is a subgroup of $M_\phi$ by Corollary \ref{subgroup}, then
\[ [M_\phi:G]=[M_\phi:M_\psi][M_\psi:G] \]
(Theorem 4.5, p.39 in \cite{HU}, where $[N:H]$ is the index of the subgroup $H$ in $N$).
Since $[M_\phi:G]=k_1k_2\cdot\cdot\cdot k_l\in{\mathbb Z}^+$, it follows that $[M_\phi:M_\psi]$ is finite.
\end{proof}

\begin{theorem}\label{algebraic} If $\phi$ is a quasiperiodic flow on $T^n$ that is algebraic, if $\psi$ is a flow on $T^n$ that is smoothly conjugate to a flow generated by constant vector field, and if $\phi$ is smoothly semiconjugate to $\psi$, then $\psi$ is a quasiperiodic flow that is algebraic, and $M_\psi$ is a finite index subgroup of $M_\phi$.
\end{theorem}

\begin{proof} The hypotheses and Theorem \ref{covering} imply that $\psi$ is a quasiperiodic flow that is algebraic. Since $\phi$ is algebraic, there is a real algebraic number field $F$ of degree $n$ such that $M_\phi$ is a subgroup of ${\mathfrak o}^*_F$ (Theorem 3.4 in \cite{BA3}). The group ${\mathfrak o}^*_F$ is finitely generated (Dirichlet's Unit Theorem), so that $M_\phi$ is finitely generated (Corollary 1.7, p.74 in \cite{HU}). By Corollary \ref{finite}, $M_\psi$ is a finite index subgroup of $M_\phi$.
\end{proof}

An application of Theorem \ref{algebraic} gives the existence of multipliers of a quasiperiodic flow $\psi$ other than $\pm 1$ when there is an algebraic quasiperiodic flow $\phi$ with multipliers other than $\pm 1$ such that $\phi$ is smoothly semiconjugate to $\psi$. (This sufficient condition for the existence of multipliers other than $\pm 1$ was implicitly used in Example \ref{example}.)

\begin{corollary} Suppose that $\phi$ is a quasiperiodic flow on $T^n$ that is algebraic, with $X_\phi$ a constant vector field and $M_\phi\setminus\{1,-1\}\ne\emptyset$. If $\psi$ is the flow on $T^n$ generated by $X_\psi=BX_\phi$ for any $B\in{\rm M}(n,{\mathbb Z})\cap{\rm GL}(n,{\mathbb Q})$, then $M_\psi\setminus\{1,-1\}\ne\emptyset$.
\end{corollary}

\begin{proof} Any $B\in{\rm M}(n,{\mathbb Z})\cap{\rm GL}(n,{\mathbb Q})$ induces a surjective $V\in C^\infty(T^n)$, so that if $X_\psi=BX_\phi={\bf T}V X_\phi$, then $[M_\phi:M_\psi]$ is finite by Theorem \ref{algebraic}. If $M_\phi\setminus\{1,-1\}\ne\emptyset$, then $M_\phi$ has an infinite cyclic factor, which implies that $M_\psi\setminus\{1,-1\}\ne\emptyset$.
\end{proof}

\end{document}